\documentclass[10t]{article}
\usepackage{amssymb}
\usepackage{mathrsfs}
\usepackage{amsfonts}
\usepackage{tikz}
\usepackage[numbers,sort&compress]{natbib}

\usepackage{amsmath,amsthm,amsfonts,amssymb,amscd}

\allowdisplaybreaks[4]
\newtheorem {Problem} {Problem}[section]
\newtheorem {Theorem} [Problem]{Theorem}
\newtheorem {Lemma}[Problem]{Lemma}

\newenvironment {Proof}{\noindent {\bf Proof.}}{\hfill\ensuremath{\square}}
\newcommand*{\QEDB}{\hfill\ensuremath{\square}}

\begin{document}

\title{The signless Laplacian spectral radius of graphs  with forbidding linear forests   \thanks{This work is supported by the National Natural Science Foundation of China (Nos. 11971311 and 11531001) and  the Montenegrin-Chinese
Science and Technology Cooperation Project (No.3-12).
\newline \indent Email: chenmingzhuabc@163.com, aming8809@163.com, xiaodong@sjtu.edu.cn. \newline  \indent $^{\dagger}$Corresponding author:
Xiao-Dong Zhang (Email: xiaodong@sjtu.edu.cn),}}

\author{ Ming-Zhu Chen, A-Ming Liu,\\
School of Science, Hainan University, Haikou 570228, P. R. China, \\
\and  Xiao-Dong Zhang$^{\dagger}$
\\School of Mathematical Sciences, MOE-LSC, SHL-MAC\\
Shanghai Jiao Tong University,
Shanghai 200240, P. R. China}

\date{}
\maketitle

\begin{abstract}
 Tur\'{a}n type extremal problem  is how to  maximize the number of edges over all graphs which do not
contain fixed forbidden subgraphs. Similarly, spectral Tur\'{a}n type extremal problem  is how to  maximize (signless Laplacian)  spectral radius over all graphs which do not contain fixed  subgraphs.
In this paper,  we first present a stability result for $k\cdot P_3$ in terms of the number of edges  and then  determine all extremal graphs maximizing the signless Laplacian spectral radius over all graphs which do not contain a fixed linear forest with at most two odd paths or  $k\cdot P_3$ as a subgraph, respectively.
\end{abstract}

\section{Introduction}
 Let $G$ be an undirected simple graph with vertex set
$V(G)=\{v_1,\dots,v_n\}$ and edge set $E(G)$, where $e(G)$ is the number of edges of $G$.
The \emph{adjacency matrix}
of $G$  is the $n\times n$ matrix $A(G)=(a_{ij})$, where
$a_{ij}=1$ if $v_i$ is adjacent to $v_j$, and $0$ otherwise.
The matrix $Q(G)=D(G)+A(G)$ is known as  the \emph {signless Laplacian matrix} of $G$,  where $D(G)$ is the degree diagonal matrix of $G$.  The \emph{spectral radius} and \emph{signless Laplacian spectral radius} of $G$ are the largest eigenvalues of $A(G)$ and $Q(G)$, denoted by
$\rho(G)$ and $q(G)$, respectively.
For $v\in V(G)$,  the \emph{neighborhood} $N_G(v)$ of $v$  is $\{u: uv\in E(G)\}$ and the \emph{degree} $d_G(v)$ of $v$  is $|N_G(v)|$.
We write $N(v)$ and $d(v)$ for $N_G(v)$ and $d_G(v)$ respectively if there is no ambiguity.
 Denote  by $\triangle(G)$ and $\delta(G)$  the maximum and minimum degree of $G$, respectively.
For $V_1, V_2\subseteq V(G)$,  $e(V_1,V_2)$ denotes the number of the  edges of  $G$ with one end vertex in $V_1$ and the other  in $V_2$. We say a graph $G$ is \emph{$F$-free} if   it  does not contain $F$ as a subgraph.
A path of order $n$ is denoted by $P_n$. A \emph{linear forest} is a forest whose
components are paths. For a path $P_3$, say $xyz$, we call $y$ its center and $x,z$ its two ends.
For two vertex disjoint graphs $G$ and $H$,  we denote by  $G\cup H$ and  $G\nabla H$  the \emph{union} of $G$ and $H$,
and the \emph{join} of $G$ and $H$, i.e., joining every vertex of $G$ to every vertex of $H$, respectively.
Denote by $k\cdot G$  the union of  $k$ disjoint  copies of $G$.
For graph notation and terminology undefined here, we refer the readers to \cite{BM}.

%
The problem of maximizing the number of edges over all graphs without  fixed subgraphs is one of the cornerstones of graph theory. In 2010,  Nikiforov \cite{Nikiforov} proposed the following
 spectral extremal   problem,  which is the spectral analogue of Tur\'{a}n type extremal problem.

\begin{Problem}\cite{Nikiforov}\label{P1}
Given a graph $H$, what is the maximum $\rho(G)$ of a graph $G$ of order $n$ which does not contain $H$ as a subgraph$?$
\end{Problem}


A variation of Problem \ref{P1} in terms of signless Laplacian spectral radius \cite{NY}  is as follows:

\begin{Problem}\label{P2}
 Given a graph $H$, what is the maximum $q(G)$ of a graph  $G$  of order $n$ which does not contain $H$ as a subgraph$?$
\end{Problem}

Recently, the signless Laplacian spectral radius of graphs without fixed subgraphs has received more and more attention. For example,  He, Jin, and Zhang \cite{JZ} obtained the signless Laplacian spectral analogue of  Tur\'{a}n  theorem (see \cite{Turan}).   Nikiforov and Yuan \cite{NY} obtained the signless Laplacian spectral analogue of Erd\H{o}s--Gallai theorem (see \cite{EG}).
For more details,  readers may be referred to \cite{CZ2,FNP,FNP2,Nikiforov2010,Nikiforov2011,NY2,Tait2017,Tait2019,Yuan}.

A natural extension of the signless Laplacian spectral analogue of
Erd\H{o}s--Gallai theorem is to determine the maximum signless Laplacian spectral radius
of graphs which do not contain a fixed linear forest as a subgraph. 
 Lidick\'{y}, Liu, and Palmer \cite{LLP}  determined the Tur\'{a}n number for a forbidden linear forest except for $ k\cdot P_3$  if the order of graph is  sufficiently large  and all extremal graphs.  Bushaw and Kettle \cite{BK}, Campos and Lopes \cite{campos2015}, and Yuan and Zhang \cite{YZ}, independently,  determined  the Tur\'{a}n number for a forbidden $ k\cdot P_3$.
 In order to state these results, we need some symbols for given graphs.

Let $S_{n,h}$ be the graph obtained from  $K_h\cup \overline{K}_{n-h}$ by adding all edges between $K_h$ and  $\overline{K}_{n-h}$, i.e., $S_{n,h}=K_h\nabla \overline{K}_{n-h}$. Moreover,
let $S^+_{n,h}$ be the graph obtained by adding an edge to $S_{n,h}$ (see Fig.~\ref{fig1}) .

\begin{center}\label{fig1}
\begin{tikzpicture}[scale=1.2]
\draw (0,0) ellipse (1cm and 0.5cm);\draw (0,-1.5) ellipse (1cm and 0.5cm);
\path (-0.6,-0.2)  coordinate (P1);\path (0.6,-0.2)  coordinate (P2);
\path (-0.6,-1.3)  coordinate (P3);\path (0.6,-1.3)  coordinate (P4);
\path (-0.2,-0.2)  coordinate (Q1);\path (0,-0.2)  coordinate (Q2);\path (0.2,-0.2)  coordinate (Q3);
\path (-0.2,-1.3)  coordinate (Q4);\path (0,-1.3)  coordinate (Q5);\path (0.2,-1.3)  coordinate (Q6);
\foreach \i in {1,2,3,4}
{\fill (P\i) circle (1pt);}
\foreach \i in {1,...,6}
{\fill (Q\i) circle (0.6pt);}
\foreach \i in {1,2}
{
\foreach \j in {3,4}
{
\draw (P\i) -- (P\j);
}
}
\path (0,0.2) node () {\scriptsize{$\,K_h$\,}};\path (0,-1.7) node () {\scriptsize{$\,\overline{K}_{n-h}$\,}};
\path (0,-2.4) node () {$\,S_{n,h}$\,};

\draw (3,0) ellipse (1cm and 0.5cm);\draw (3,-1.5) ellipse (1cm and 0.5cm);
\path (2.4,-0.2)  coordinate (P1);\path (3.6,-0.2)  coordinate (P2);
\path (2.4,-1.3)  coordinate (P3);\path (3.6,-1.3)  coordinate (P4);
\path (2.8,-1.3)  coordinate (P5); \path (3,-1.3)  coordinate (P6);
\path (2.8,-0.2)  coordinate (Q1);\path (3,-0.2)  coordinate (Q2);\path (3.2,-0.2)  coordinate (Q3);
\path (3.2,-1.3)  coordinate (Q4);\path (3.3,-1.3)  coordinate (Q5);\path (3.4,-1.3)  coordinate (Q6);
\foreach \i in {1,...,6}
{\fill (P\i) circle (1pt);}
\foreach \i in {1,...,6}
{\fill (Q\i) circle (0.6pt);}
\foreach \i in {1,2}
{
\foreach \j in {3,...,6}
{
\draw (P\i) -- (P\j);
}
}
\draw (P3) -- (P5);
\path (3,0.2) node () {\scriptsize{$\,K_h$\,}};\path (3,-1.7) node () {\scriptsize{$\,K_2\cup\overline{ K}_{n-h-2}$\,}};
\path (3,-2.4) node () {$\,S^+_{n,h}$\,};

\path (2,-3) node () { Fig. 1. $S_{n,h}$ and $S_{n,h}^+$};
\end{tikzpicture}
\end{center}

  Let $L_{t_1,t_2,h,h+1} $ be the graph obtained from $t_1K_{h}\cup t_2K_{h+1}$ by joining all edges from one new vertex to all vertices of $t_1K_{h}\cup t_2K_{h+1},$  
 where $t_1,t_2\geq0$,  i.e.,
   $L_{t_1,t_2,h,h+1}=K_1\nabla(t_1K_h\cup t_2K_{h+1})$.
In particular, write $L_{t,h}=L_{t,0,h,h+1}$ (see Fig.~2).

\begin{center}
\label{fig2}
\begin{tikzpicture}[scale=1.2]
\path (1,0)  coordinate (P1);\path (0,0.3)  coordinate (P2);\path (0.5,1)  coordinate (P3);
\path (1.5,1)  coordinate (P4);\path (2,0.3)  coordinate (P5);
\path (0,-0.3)  coordinate (P6);\path (0.5,-1)  coordinate (P7);
\path (1.5,-1)  coordinate (P8);\path (2,-0.3)  coordinate (P9);
\path (0.8,0.7)  coordinate (Q1);\path (1,0.7)  coordinate (Q2);\path (1.2,0.7)  coordinate (Q3);
\path (0.8,-0.7)  coordinate (Q4);\path (1,-0.7)  coordinate (Q5);\path (1.2,-0.7)  coordinate (Q6);
\foreach \i in {1,...,9}
{\fill (P\i) circle (1pt);}
\foreach \i in {1,...,6}
{\fill (Q\i) circle (0.5pt);}
\draw (P1)--(P2)--(P3)--cycle;
\draw (P1)--(P4)--(P5)--cycle;
\draw (P1)--(P6)--(P7)--cycle;
\draw (P1)--(P8)--(P9)--cycle;
\path (0.5,0.4) node () {\scriptsize{$\,K_{h+1}$\,}};\path (1.6,0.4) node () {\scriptsize{$\,K_{h+1}$\,}};
\path (0.5,-0.4) node () {\scriptsize{$\,K_{h+2}$\,}};\path (1.6,-0.4) node () {\scriptsize{$\,K_{h+2}$\,}};
\path (1,0.9) node () {$\,t_1$\,};\path (1,-0.9) node () {$\,t_2$\,};
\path (1,-1.5) node () {$\,L_{t_1,t_2,h,h+1}$\,};

\path (5,-0.5)  coordinate (P1);\path (4,-0.2)  coordinate (P2);\path (4.5,0.5)  coordinate (P3);
\path (5.5,0.5)  coordinate (P4);\path (6,-0.2)  coordinate (P5);
\path (4.8,0.2)  coordinate (Q1);\path (5,0.2)  coordinate (Q2);\path (5.2,0.2)  coordinate (Q3);
\foreach \i in {1,...,5}
{\fill (P\i) circle (1pt);}
\foreach \i in {1,...,3}
{\fill (Q\i) circle (0.5pt);}
\draw (P1)--(P2)--(P3)--cycle;
\draw (P1)--(P4)--(P5)--cycle;
\path (4.5,-0.1) node () {\scriptsize{$\,K_{h+1}$\,}};\path (5.6,-0.1) node () {\scriptsize{$\,K_{h+1}$\,}};
\path (5,0.4) node () {$\,t$\,};
\path (5,-1.5) node () {$\,L_{t,h}$\,};
\path (3,-2.3) node () {Fig. 2. $L_{t_1,t_2,h,h+1}$ and $L_{t,h}$};
\end{tikzpicture}
\end{center}

Suppose that $H$ is a  connected graph  of order $h$. For $k\geq1$, $h\geq2$, and $n>k+h-1$,  we define
$F_{n,k}(H):=K_{k-1}\nabla(H\cup p\cdot K_2\cup K_s)$, where $n-(k+h-1)=2p+s$ and $0\leq s<2$.
In particular,  if $H=K_2$, then we write $F_{n,k}$ for $F_{n,k}(K_2)$.
Let $H_{n,1}$ be a graph of order  $n\geq 7$ obtained from $S_{n-2,2}$ and $K_3$ by identifying a vertex of maximum degree in $S_{n-2,2}$ with
a vertex of $K_3$.
Let $N_6$ be a graph of order $6$ by adding a pendant edge to every vertex of $K_3$, respectively (see Fig.~3).

\begin{center}
\begin{tikzpicture}[scale=1.4]
\path (2.5,0)  coordinate (P1);\path (4,0)  coordinate (P2);
\path (3.25,0.4)  coordinate (P3);\path (3.25,0.75)  coordinate (P4);\path (3.25,-0.75)  coordinate (P5);
\path (3.25,-0.3)  coordinate (Q1);\path (3.25,-0.4)  coordinate (Q2);\path (3.25,-0.5)  coordinate (Q3);
\path (4.5,0.5)  coordinate (R1);\path (4.5,-0.5)  coordinate (R2);
\foreach \i in {1,2,3,4,5}
{\fill (P\i) circle (1pt);}
\foreach \i in {1,2,3}
{\fill (Q\i) circle (0.7pt);}
\fill (R1) circle (1pt);\fill (R2) circle (1pt);
\foreach \i in {1,2}
{
\foreach \j in {3,4,5}
{
\draw (P\i) -- (P\j);
}
}
\draw (P1) -- (P2);
\draw (P2) -- (R1)--(R2)--cycle;
\path (3.3,-1.3) node () {$H_{n,1}$};

\draw (7,0.5) ellipse (1cm and 0.2cm);\draw (7,-0.6) ellipse (1cm and 0.2cm);\draw (7.6,-0.6) ellipse (0.1cm and 0.1cm);
\path (6.4,0.5)  coordinate (P1);\path (7.6,0.5)  coordinate (P2);
\path (6.2,-0.6)  coordinate (P3);\path (7.6,-0.6)  coordinate (P4);
\path (6.6,-0.6)  coordinate (P5); \path (6.8,-0.6)  coordinate (P6);\path (7,-0.6)  coordinate (P7);
\path (6.8,0.5)  coordinate (Q1);\path (7,0.5)  coordinate (Q2);\path (7.2,0.5)  coordinate (Q3);
\path (7.2,-0.6)  coordinate (Q4);\path (7.3,-0.6)  coordinate (Q5);\path (7.4,-0.6)  coordinate (Q6);
\foreach \i in {1,...,7}
{\fill (P\i) circle (1pt);}
\foreach \i in {1,...,6}
{\fill (Q\i) circle (0.6pt);}
\foreach \i in {1,2}
{
\foreach \j in {3,...,7}
{
\draw (P\i) -- (P\j);
}
}
\draw (P3) -- (P5);\draw (P6) -- (P7);
\path (7,0.85) node () {\scriptsize{$K_{k-1}$}}; \path (7.8,-0.8) node () {\scriptsize{$K_{s}$}};
\path (7,-1.3) node () {$F_{n,k}$};

\path (9.4,-0.6)  coordinate (P1);\path (10.8,-0.6)  coordinate (P2); \path (10.1,0.2)  coordinate (P3);
\path (8.9,-0.6)  coordinate (Q1);\path (11.3,-0.6)  coordinate (Q2); \path (10.1,0.7)  coordinate (Q3);
\fill (P1) circle (1pt); \fill (P2) circle (1pt);\fill (P3) circle (1pt);
\fill (Q1) circle (1pt); \fill (Q2) circle (1pt);\fill (Q3) circle (1pt);
\draw (P1)--(P2)--(P3)--cycle;
\draw (P1)--(Q1);\draw (P2)--(Q2);\draw (P3)--(Q3);
\path (10.2,-1.3) node () {$N_6$};

\path (6.5,-2) node () {Fig. 3. $H_{n,1}$, $F_{n,k}$, and $N_6$};
\end{tikzpicture}
\end{center}


\begin{Theorem}\label{edge extremal graph for linear forest}\cite{LLP}
Let $F=\cup_{i=1}^k P_{a_i}$ be a linear forest  with $k\ge 2 $, $a_1\geq \cdots \geq a_k\geq2$,  and  $h=\sum_{i=1}^k \lfloor\frac{a_i}{2}\rfloor-1$. If   there exists at least one $a_i$ not $3$ and $G$ is  an $F$-free  graph  of  order $n$,  then for  sufficiently large $n$,
$$e(G)\leq  \binom{h}{2}+ h(n-h)+c,$$
where $c=1$ if all $a_i$ are odd and $c=0$ otherwise. Moreover,
if $c=1$ then  the equality holds if and only if $G=S^+_{n,h}$. Otherwise, the  equality holds if and only if  $G=S_{n,h}$.
\end{Theorem}

\begin{Theorem}\label{edge extremal graph for kP3}\cite{BK, campos2015,YZ}
Let  $G$ be a \ $k\cdot P_3$-free graph of order $n$. Then 
\[e(G)\leq\left\{
\begin{array}{llll}
 \vspace{1mm}
 \binom{n}{2},&& \mbox{for $n<3k$};\\
   \vspace{1mm}
  \binom{3k-1}{2}+\big\lfloor\frac{n-3k+1}{2}\big\rfloor,&& \mbox{for $3k\leq n<5k-1$};\\
       \vspace{1mm}
 \binom{3k-1}{2}+k, && \mbox{for $n=5k-1$};\\
      \vspace{1mm}
 \binom{k-1}{2}+(n-k+1)(k-1)+\big\lfloor\frac{n-k+1}{2}\big\rfloor, && \mbox{for $n>5k-1$.}
\end{array}\right.
\]
Moreover, (i) If $n <3k$, then the equality holds if and only if $G=K_n$;\\
(ii) If $3k\leq n < 5k-1$, then the equality holds if and only if  $G=K_{3k-1}\cup F_{n-3k+1,1}$;\\
(iii) If $n = 5k - 1$, then the equality holds if and only if  $G=K_{3k-1}\cup F_{2k,1}$ or $G=F_{5k-1,k}$;\\
(iv) If $n > 5k - 1$, then the equality holds if and only if $G=F_{n,k}$.
\end{Theorem}


In this paper, we have made some contributions to Tur\'{a}n type results and Problem~\ref{P2}.  Firstly, we obtain a stability result for $k\cdot P_3$ in terms of the number of edges, which improves  Theorem~\ref{edge extremal graph for kP3}.  Then we   determine all extremal graphs maximizing the signless Laplacian spectral radius over  all  graphs of  sufficiently large order $n$ without  a  fixed linear forest with at most two odd paths  or $ k\cdot P_3$ as a subgraph.    Our main results are stated as follows:

\begin{Theorem}\label{stablility thm for kP3}
Let $G$ be a graph of order $n\geq\frac{11}{2}k^2+2k-\frac{3}{2}$ with $k\geq2$. If $e(G)>(k-\frac{3}{2})n$, then $G$ contains $k\cdot P_3$ as a subgraph unless one of the following holds: \\
(i) $G\subseteq F_{n,k}$;\\
(ii) $G\subseteq F_{n,k-1}(K_h)$, where $4\leq h\leq5$;\\
(iii) $G\subseteq F_{n,k-1}(N_6)$.
\end{Theorem}

\noindent{\bf Remark.} Roughly speaking, by Theorem~\ref{edge extremal graph for kP3}, it is easy to see that  a graph of sufficiently large order  $n$ with at least $(k-\frac{1}{2})n$ edges contains $k\cdot P_3$ as a subgraph. But, by Theorem~\ref{stablility thm for kP3}, a graph of order  $n$ with at least $(k-\frac{3}{2})n$ edges contains $k\cdot P_3$  except several known graphs as a subgraph. Hence  Theorem~\ref{stablility thm for kP3} improves  Theorem~\ref{edge extremal graph for kP3} and the extremal graph $F_{n,k}$ for $k\cdot P_3$ is stable.

 \begin{Theorem}\label{Q spectral radius for all linear forest}
 Let  $F=\cup_{i=1}^k P_{a_i}$ be a linear forest with   $k\geq2$ and
$a_1\geq \cdots \geq a_k\geq2$ and $ h=\sum _{i=1}^k \lfloor\frac{a_i}{2}\rfloor-1$.  Suppose that $G$ is  an $F$-free  graph  of sufficiently large order $n$. \\
(i) If there  exists  at least one even $a_i$ and  at most  two odd  $a_j$, then $q(G)\leq q(S_{n,h})$ with equality if and only if $G=S_{n,h}$;\\
(ii) If  $k=2$ and  both $a_1>3$ and $a_2\geq3$ are  odd,  then $q(G)\leq q(S^+_{n,h})$ with equality if and only if $G=S^+_{n,h}$.\\
(iii) If $a_1=\cdots=a_k= 3,$ i.e., $F=k\cdot P_3$, then $q(G)\leq q(F_{n,k}) $ with equality if and
only if   $G=F_{n,k}$.
\end{Theorem}


The rest of this paper is organized as follows. In Section~2, some known and new results are presented.
In Section~3, we give  the proof of Theorem~\ref{stablility thm for kP3}.
In Section~4, we give  the proof of Theorem~\ref{Q spectral radius for all linear forest}.

\section{Preliminary}

We first present some stability results of graphs with large minimum degree and without  a fixed linear forest as a subgraph,  which play an important role in the proof of our main results.

\begin{Theorem}\label{Thm2.1}\cite{CZ}
Let $F=\cup_{i=1}^k P_{a_i}$ be a linear forest with   $k\geq2$ and
$a_1\geq \cdots \geq a_k\geq2$. Denote $ h=\sum_{i=1}^k \lfloor\frac{a_i}{2}\rfloor-1$ and suppose that $G$ is  an $F$-free connected graph  of order $n\geq2k+2$.
 If all $a_i$ are even and  $\delta(G)\geq h$,
 then  one of the following holds:\\
 (i) $G\subseteq S_{n,h}$;\\
  (ii) $F=2\cdot P_{a_1}$ and $G=L_{t,h}$, where $n=th+1$.
\end{Theorem}

\begin{Theorem}\label{Thm2.2}\cite{CZ}
Let  $F=\cup_{i=1}^k P_{a_i}$ be a linear forest with    $k\geq2$ and
$a_1\geq \cdots \geq a_k\geq2$.  Denote $ h=\sum_{i=1}^k \lfloor\frac{a_i}{2}\rfloor-1$ and suppose that $G$ is  an $F$-free connected graph  of order $n\geq2k+3$.
 If there exists precisely one odd  $a_i$  and  $\delta(G)\geq h$,
then  one of the following holds:\\
 (i) $G\subseteq S_{n,h}$;\\
 (ii)  $F=P_6\cup P_3$ and $G\subseteq F_{n,3}$, where $n$ is even;\\
  (iii)  $F\in \{P_{a_{1}-1}\cup P_{a_{1}}, P_{a_{1}+1}\cup P_{a_{1}}\}$ and $G= L_{t,h}$, where $a_1$ is odd and  $n=th+1$.
\end{Theorem}

\begin{Theorem}\label{Thm2.3}\cite{CZ}
Let  $F=\cup_{i=1}^k P_{a_i}$ be a linear forest with   $k\geq2$ and
$a_1\geq \cdots \geq a_k\geq2$.  Denote $ h=\sum_{i=1}^k \lfloor\frac{a_i}{2}\rfloor-1$ and suppose that $G$ is  an $F$-free $2$-connected graph  of order  $n\geq4(2h+1)^2\binom{2h+1}{h}$. Suppose that there exist precisely two odd $a_i$.

(a) If  $\delta(G)\geq h$ and $k=2$, then  one of the following holds:\\
 (i) $G\subseteq S^+_{n,h}$;\\
  (ii) $F=P_7\cup P_3$ and $G\subseteq F_{n,3}$, where $n$ is even;\\
   (iii)  $F=P_9\cup P_3$ and $G\subseteq F_{n,4}$, where $n$ is odd.

(b) If  $\delta(G)\geq h$ and  $k\geq3$, then one of the following holds:\\
(iv) $G\subseteq S_{n,h}$;\\
  (v) $F= P_4\cup 2\cdot P_3$ and $G\subseteq F_{n,3}$, where $n$ is even;\\
   (vi)  $F=P_6\cup 2\cdot P_3$ and $G\subseteq F_{n,4}$, where  $n$ is odd.
\end{Theorem}

\begin{Theorem}\label{Thm2.4}\cite{CZ}
Let  $F=\cup_{i=1}^k P_{a_i}$ be a linear forest with    $k\geq2$ and
$a_1\geq \cdots \geq a_k\geq2$.  Denote $ h=\sum_{i=1}^k \lfloor\frac{a_i}{2}\rfloor-1$ and suppose that $G$ is  an $F$-free connected graph  of order $n\geq2k+4$ with at least one cut vertex. Suppose that there exist precisely two odd  $a_i$  and  $\Delta(G)=n-1$.

(a)  If   $\delta(G)\geq h\geq1$ and  $k=2$, then  one of the following holds:\\
 (i) $F=P_5\cup P_3$ and  $G\subseteq H_{n,1}$;\\
 (ii) $F=P_{a_{1}}\cup P_{a_{1}-2}$ and $G=L_{t,h}$, where  $a_1$ is odd and $n=th+1$;\\
 (iii) $F=2\cdot P_{a_{1}}$ and $G\subseteq L_{t_1,t_2,h,h+1}$, where  $a_1$ is odd and  $n=t_1h+t_2(h+1)+1$.

 (b) If $\delta(G)\geq h\geq2$ an $k\geq3$, then   $F=P_2\cup 2\cdot P_{a_{1}}$ and $G=L_{t,h}$, where $a_1$ is odd and $n=th+1$.
\end{Theorem}

\begin{Theorem}\label{Thm2.6}\cite{EG} (Erd\H{o}s--Gallai theorem)
Let $l\geq2$ and $G$ be a   graph of order $n$.   If $e(G)>\frac{(l-2)n}{2}$,  then $G$ contains a path of order $l$  as a subgraph.
\end{Theorem}

Next we give several  lemmas regarding signless Laplacian spectral radius.

\begin{Lemma}\label{Lem2.2}\cite{NY}
Let $ h\geq2$ and  $n\geq 7h^2$. Then \\
(i)  $q(S^+_{n,h})>q(S_{n,h}) >n+2h-2-\frac{2(h^2-h)}{n+2h-3}>n+2h-3$.\\
  (ii)    $q(L_{t,h})<q(S_{n,h})$, where $n=th+1$.\\
  (iii) If $G$ is a graph of order $n$ with  $q(G)\geq q(S_{n,h})$, then  $e(G)>h(n-h)$.
\end{Lemma}

\begin{Lemma}\label{Lem2.3}
(i) For $n\geq28$, $q(H_{n,1})<q(S_{n,2})$.\\
(ii) For  $h\geq3$, $t_1,t_2\geq 0$, and $n=t_1h+t_2(h+1)+1\geq 7h^2$, $q(L_{t_1,t_2,h,h+1})<q(S_{n,h})$.
\end{Lemma}

\begin{Proof}
(i) Denote $q=q(H_{n,1})$ and $Q=Q(H_{n,1})$.
Let $\mathbf x$ be a positive  eigenvector of $Q$ corresponding to $q$. Denote by $x_1$ and $x_2$ the eigenvector entry of the vertex $u$ with the maximum  degree $n-1$ and  the vertex $v$ with the the second largest  degree $n-3$, respectively. By symmetry, those two vertices $w_1, w_2$ with degree $2$ not adjacent to $v$ have the same eigenvector entry, denoted by $x_3$. Similarly, by symmetry, the remaining vertices in $H_{n,1}$ also have the same eigenvector entry, denoted by $x_4$.

Since $S_{n-2,2}$ is a proper subgraph of $H_{n,1}$, by Perron--Frobenius theorem and Lemma~\ref{Lem2.2}, $q>q(S_{n-2,2})>n-1$. By eigenequation of $Q$, we have
\[
\begin{array}{llll}
 \vspace{1mm}
 (q-n+1)x_1&=& x_2+2x_3+(n-4)x_4,\\
   \vspace{1mm}
  (q-n+3)x_2&=&x_1+(n-4)x_4,\\
       \vspace{1mm}
 (q-2)x_3&=&x_3+x_1,\\
      \vspace{1mm}
(q-2)x_4&=&x_1+x_2.
\end{array}
\]
By simple calculation, we have
$$\big((q-n+2)(q-3)-2\big)x_3=(q-n+4)x_2,$$
which implies that
$$\frac{x_2}{x_3}=\frac{(q-n+2)(q-3)-2}{q-n+4}\geq1$$
as $q>n-1$ and $n\geq 9$. Hence $x_2\geq x_3$.

Let $G$ be a graph obtained from  $H_{n,1}$ by deleting the edge $ w_1w_2$ and adding the edges
$\{vw_1,vw_2\}$.
Since\begin{eqnarray*}
       q(G)-q  &\geq&   \mathbf x^{T}Q(G)\mathbf x- \mathbf x^T Q\mathbf x\\
       &=& 2(x_2+x_3)^2-(x_3+x_3)^2\\
       &=&2(x_2-x_3)(x_2+x_3)+4x_2x_3 >0,
     \end{eqnarray*}
     we have $ q(G)>q$.  Note that $G=S_{n,2}$. Hence $ q<q(S_{n,2})$.

(ii) Note that $q(G)\leq \max \limits_{uv\in E(G)}  \Big\{d(u)+d(v) \Big\}$ \cite{DPS}. By Lemma~\ref{Lem2.2}~($i$),
\begin{eqnarray*}
  q(L_{t_1,t_2,h,h+1}) \leq \max \limits_{uv\in E(L_{t_1,t_2,h,h+1})} \Big\{d(u)+d(v)\Big\}
   \leq n+h
   \leq n+2h-3
  < q(S_{n,h}).
\end{eqnarray*}

\end{Proof}

\begin{Lemma}\label{lbound for extremal graph}
For $k\geq2$ and $n\geq 2k^2$, $$n+2k-5<q(F_{n,k})\leq \frac{n+2k-2+\sqrt{(n+2k-6)^2-8(k^2-4k+3)}}{2}.$$
Moreover, if $k\geq3$ then $q(F_{n,k})\leq n+2k-4$.
\end{Lemma}


\begin{Proof}
Denote $q=q(F_{n,k})$ and $Q=Q(F_{n,k})$.
 Let $\mathbf x$ be a positive eigenvector of $Q$ corresponding to $q$.
Let $n-(k-1)=2p+s$ with $0\le s<2$. By symmetry, all vertices of  subgraphs $K_{k-1}$, $p\cdot K_2$, or $K_s$ in  $F_{n,k}=K_{k-1}\nabla(p\cdot K_2 \cup K_s)$ have the same  eigenvector entry respectively, which are denoted by $x_1$, $x_2$, $x_3$, respectively.

(i) If $n-k+1$ is even, then $s=0$.
By  $Q\mathbf x=q\mathbf x$, it is easy to see that
\begin{eqnarray*}
    (q-n+1) x_1 &=& (k-2)x_1 +(n-k+1)x_2, \\
    (q-k) x_2 &=& (k-1)x_1+x_2.
\end{eqnarray*}
It is easy to see that $q$ is the largest root of
 $f(x)=0$, where $f(x)=x^2-(n+2k-2)x+2n+2k^2-4k-2=0$.
Hence $$q=\frac{n+2k-2+\sqrt{(n+2k-6)^2-8(k^2-4k+3)}}{2}>n+2k-5.$$
Moreover, if $k\geq3$ then $q\leq n+2k-4$.

(ii) If $n-k+1$ is odd, then $s=1$. By  $Q\mathbf x=q\mathbf x$, it is easy to see that
\begin{eqnarray*}
    (q-n+1)x_1 &=& (k-2)x_1 +(n-k)x_2+x_3, \\
     (q-k)x_2 &=& (k-1)x_1+x_2,    \\
     (q-k+1)x_3 &=& (k-1)x_1.
\end{eqnarray*}
It is easy to see that $q$ is the largest root of $f(x)=0$, where 
$$f(x)=(x-k+1)[x^2-(n+2k-2)x+2n+2k^2-4k-2]+2k-2.$$
 Note that  $q>q(K_k)=2k-2$,
 we have
 $$q^2-(n+2k-2)q+2n+2k^2-4k-2=-\frac{2k-2}{q-k+1}<0,$$
 which implies that
$$q<\frac{n+2k-2+\sqrt{(n+2k-6)^2-8(k^2-4k+3)}}{2}.$$
Furthermore, if $k\geq3$ then $q\leq n+2k-4$.
Moreover,
since $$f(n+2k-5)=-n^2+(2k^2-11k+17)n+2k^3-18k^2+55k-54<0,$$ we have $q>n+2k-5$.
\end{Proof}

\begin{Lemma}\label{Kh}
For $k\geq2$, $4\leq h\leq5$ and  $n\geq k+20$, $q(F_{n,k-1}(K_h))<q(F_{n,k})$.
\end{Lemma}

\begin{Proof}
Denote $q=q(F_{n,k-1}(K_4))$ and $Q=Q(F_{n,k-1}(K_4))$. Let $\mathbf x$ be a positive  eigenvector of $Q$ corresponding to $q$.
Let $V(K_4)=\{v_1,\dots,v_4\}$  and $V(p\cdot K_2)=\{w_1,\dots,w_{2p}\}$  in the representation of $F_{n,k-1}(K_4)$, where $p=\lfloor\frac{n-k-2}{2}\rfloor$  and $w_{2i-1}$ is adjacent to $w_{2i}$ for $1\leq i\leq p$.
By symmetry, all vertices of $K_4$ and $p\cdot K_2$ in the representation of $F_{n,k-1}(K_4)$ have the same eigenvector entry respectively, which are denoted by $x_1$ and $x_2$ respectively.
Let $G$ be a graph obtained from  $F_{n,k-1}(K_4)$ by deleting edges $ \{v_2v_4,v_3v_4\}$ and adding all the edges $v_1w_j$ for $1\leq j\leq2\big\lfloor\frac{n-k-2}{2}\big\rfloor$.
Since
\begin{eqnarray*}
  q(G)-q &\geq &  \mathbf x^{T}Q(G)\mathbf x- \mathbf x^T Q\mathbf x\geq (n-k-3)(x_1+x_2)^2-2(x_1+x_1)^2 \\
 &\geq& 8(x_1+x_2)^2-2(x_1+x_1)^2=8x_2(2x_1+x_2)>0,
\end{eqnarray*}
we have $ q(G)>q$.  Note that $G\subseteq F_{n,k}$, by Perron--Frobenius theorem, we have $ q<q(G)\leq q(F_{n,k})$.

Similarly, we also have $ q(F_{n,k-1}(K_5))< q(F_{n,k})$.
\end{Proof}

\begin{Lemma}\label{N1}
For $k\geq2$ and $n\geq k+4$, $q(F_{n,k-1}(N_6))<q(F_{n,k})$.
\end{Lemma}

\begin{Proof}
Denote $q=q(F_{n,k-1}(N_6))$ and $Q=Q(F_{n,k-1}(N_6))$.
Let $\mathbf x$ be a positive  eigenvector of $Q$ corresponding to $q$.
Let $V(N_6)=\{v_1,v_2,v_3,w_1,w_2,w_3\}$, where $v_1,v_2,v_3$ are the vertices of triangle in $N_6$ and $w_i$ is adjacent to $v_i$ for $1\leq i\leq 3$. By symmetry,  $x_{v_1}=x_{v_2}=x_{v_3}$, denoted by $x_1$, and
$x_{w_1}=x_{w_2}=x_{w_3}$, denoted by $x_2$.
Let $G$ be a graph obtained from  $F_{n,k-1}(N_6)$ by deleting the edges $ \{v_2w_2,v_3w_3\}$ and adding the edges
$\{v_1w_2,v_1w_3,w_2w_3\}$.
Since\begin{eqnarray*}
       q(G)-q  \geq   \mathbf x^{T}Q(G)\mathbf x- \mathbf x^T Q\mathbf x= 2(x_1+x_2)^2+(x_2+x_2)^2-2(x_1+x_2)^2
       =4x_2^2 >0,
     \end{eqnarray*}
     we have $ q(G)>q$.  Note that $G\subseteq F_{n,k}$. By the Perron--Frobenius theorem, we have $ q<q(G)\leq q(F_{n,k})$.
\end{Proof}

\section{Proof of Theorem~\ref{stablility thm for kP3}}

\begin{Lemma}\label{2P3-free}
Let $G$ be a  $2\cdot P_3$-free graph of order $n\geq6$. Then
 $G\subseteq F_{n,1}(K_h)$ with $4\leq h\leq5$,  $G\subseteq F_{n,1}(N_6)$, or  $G\subseteq F_{n,2}$.
\end{Lemma}

\begin{Proof}
If $G$ is $P_3$-free, then $G$ consists of independent edges and isolated vertices, which implies that $G\subseteq F_{n,2}$.
So we assume that  $G$ contains $P_3$ as a subgraph.
Note that every connected graph of order at least 3 contains $P_3$ as a subgraph. Since $G$ is $2\cdot P_3$-free,  $G$ must have precisely one component $H$ of order $h\geq 3$ and every of the remaining components (if any exists) is an independent edge or an isolated vertex. Obviously, $H$ is $2\cdot P_3$-free. If $3\leq h\leq 5$, then $H\subseteq K_h$,  which implies that $G\subseteq F_{n,1}(K_h)$.  So next we assume that $h\geq6$.
Let $Q=v_1v_2\cdots v_l$ be the longest path of order $l$ in $H$. Note that $ H$ contains $P_3$ as a subgraph and  $H$ is $2\cdot P_3$-free, we have $3\leq l\leq 5$.
If $l=3$, then $H$ is a star. So $G\subseteq F_{n,2}$. If $l=4$, then $V(H)\backslash V(Q)$ is an independent set of size at least 2,  in which all vertices are all adjacent to  one of $v_2$ and $v_3$, otherwise  $H$ contains $2\cdot P_3$ as a subgraph.  In addition, $v_1$ is not adjacent to $v_4$, otherwise  $H$ contains $2\cdot P_3$ as a subgraph. So $H\subseteq F_{h,2}$, which implies that  $G\subseteq F_{n,2}$.
If $l=5$, then $H-V(Q)$ consists of independent edges and isolated vertices, say $u_1u_2,\dots,u_pu_p,w_1,\dots,w_q$,
where $p+q\geq1$, in which every vertex has at most one neighbor $v_3$ in $Q$.
In addiction, since $H$ is $2\cdot P_3$-free, we have $v_1v_4, v_1v_5, v_2v_5\notin E(H)$.
If $v_2v_4\in E(H)$, then $v_1v_3, v_3v_5\notin E(H)$, $p=0$ and $q=1$,  which implies that $h=6$ and  $H=N_6$.  Thus $G\subseteq F_{n,1}(N_6)$.
If  $v_2v_4\notin E(H)$, then $H\subseteq F_{h,2}$.  Thus $G\subseteq F_{n,2}$. This completes the proof.
\end{Proof}

\vspace{3mm}
\noindent{\bf Proof of Theorem~\ref{stablility thm for kP3}.}
Suppose that $G$ is  $k\cdot P_3$-free.
Since $e(G)>(k-\frac{3}{2})n\geq e(F_{n,k-1})$,  it follows from Theorem~\ref{edge extremal graph for kP3} that $G$ contains $(k-1)\cdot P_3$ as a subgraph.
Let $U=V((k-1)\cdot P_3)$.  Since $G$ is  $k\cdot P_3$-free, $G-U$ is $P_3$ free.
Thus $G-U$ consists of independent edges and isolated vertices and $e(G-U)\leq \frac{n-3k+3}{2}$.
 Then
\begin{eqnarray*}
  e(U,V(G)\backslash U) &=& e(G)-e(U)-e(G-U) \\
   &>& (k-\frac{3}{2})n-\binom{3k-3}{2}-\frac{n-3k+3}{2}\\
   &=& (k-2)n-\frac{9k^2-24k+15}{2}.
\end{eqnarray*}
Let $C\subseteq U$ be a vertex subset such that every vertex in $C$ has at least $2k+2$ neighbors in $V(G)\backslash U$. We claim that $ k-2\leq |C|\leq k-1$. In fact, if $|C|\geq k$, then there are $k$ disjoint paths $P_3's$ with centers in $C$ and two ends in $V(G)\backslash U$ as every vertex in $C$ has at least $2k+2$ neighbors in $V(G)\backslash U$, a contradiction. If $|C|\leq k-3$, then we have
\begin{eqnarray*}
  e(U,V(G)\backslash U) &\leq& |C|(n-3k+3)+(3k-3-|C|)(2k+1) \\
   &=& |C|(n-5k+2)+(3k-3)(2k+1)\\
   &\leq&(k-3)(n-5k+2)+(3k-3)(2k+1) \\
   &=& (k-3)n+k^2+14k-9\\
   &\leq& (k-2)n-\frac{9k^2-24k+15}{2},
\end{eqnarray*}
a contradiction.

If $|C|=k-1$, then $G-C$ is $P_3$-free. Otherwise  $(k-1)\cdot P_3$ with all centers in $C$ and all ends in $V(G)\backslash U$ and another disjoint $P_3$ in $G-C$  yield $k\cdot P_3$, a contradiction.
So $G-C$ consists of independent edges and isolated vertices.  Hence $G\subseteq F_{n,k}$.

If $|C|=k-2$, then we claim that $G-C$ is $2\cdot P_3$-free. Otherwise  $(k-2)\cdot P_3$ with all centers in $C$ and all ends in $V(G)\backslash U$ and another two disjoint $P_3's$ in $G-C$  yield $k\cdot P_3$, a contradiction. Denote by $n_1$ the order of $G-C$.  By Lemma~\ref{2P3-free},  $G-C\subseteq F_{n_1,1}(K_h)$ with $4\leq h\leq5$, $G-C\subseteq F_{n_1,1}(N_6)$ or $G-C\subseteq F_{n_1,2}$, implying that $G\subseteq F_{n,k-1}(K_h)$ with $4\leq h\leq5$, $G\subseteq F_{n,k-1}(N_6)$, or $G\subseteq F_{n,k}$.
\QEDB

\section{Proof of Theorem~\ref{Q spectral radius for all linear forest}}
 In order to prove Theorem~\ref{Q spectral radius for all linear forest}, we need to prove some technical lemmas.


\begin{Lemma}\label{Lem3.1}
 Let  $F=\cup_{i=1}^k P_{a_i}$ be a linear forest  with   $k\geq2$,
$a_1\geq \cdots \geq a_k\geq2$, and  at least one  $a_i$  not equal to $3$. Denote $ h=\sum _{i=1}^k \lfloor\frac{a_i}{2}\rfloor-1$ and suppose that $G$ is  an $F$-free  graph  of sufficiently large order $n$. If $q(G)\geq q(S_{n,h})$,
then $\Delta(G)=n-1$.
\end{Lemma}

\begin{Proof}
Suppose  that $\Delta(G)\leq n-2$.  Note that
$$q(G)\leq \max \limits_{v\in V(G)}\Big\{d(v)+\frac{1}{d(v)}\sum\limits_{z\in N(v)} d(z)\Big\},$$
which dates back to Merris \cite{Merris}.
Let $u$ be a vertex such that
$$q(G) \leq \max \limits_{v\in V(G)}\Big\{d(v)+\frac{1}{d(v)}\sum\limits_{z\in N(v)} d(z)\Big\}=d(u)+\frac{1}{d(u)}\sum\limits_{z\in N(u)} d(z).$$
If $d(u)\leq2h-1$, then
$$q(G)\leq d(u)+\frac{1}{d(u)}\sum\limits_{z\in N(u)} d(z)\leq d(u)+\Delta(G)\leq n+2h-3<q(S_{n,h}),$$ a contradiction.
Hence $d(u)\geq2h$. On the other hand,
by Theorem~\ref{edge extremal graph for linear forest},
$$\sum \limits_{z\in N(u)} d(z)= 2e(G)-\sum \limits_{z\in V(G)\backslash N(u)}d(z)\leq2e(S^+_{n,h})-n+1\leq (2h-1)n-h^2-h+3.$$
Thus
$$q(G)\leq d(u)+\frac{(2h-1)n-h^2-h+3}{d(u)}.$$
Note that the function $f(x)=x+\frac{(2h-1)n-h^2-h+3}{x}$ is convex with respect to $x$ for $x>0$.
Since
$$f(2h)=n+2h-\frac{n+(h^2+h)-3}{2h}<n+2h-2-\frac{2(h^2-h)}{n+2h-3}\leq q(S_{n,h})$$
and
$$f(n-2)=n+2h-3-\frac{h^2-3h-1}{n-2}<n+2h-2-\frac{2(h^2-h)}{n+2h-3}\leq q(S_{n,h}),$$
 we have $q(G)\leq \max \Big\{f(2h),f(n-2)\Big\}<q(S_{n,h})$, a contradiction. This completes the proof.
\end{Proof}

\begin{Lemma}\label{Lem3.2}
Let  $F=\cup_{i=1}^k P_{a_i}$  be a linear forest with   $k\geq2$,
$a_1\geq \cdots \geq a_k\geq2$, and  at least one  $a_i$  not equal to $3$. Denote $ h=\sum _{i=1}^k \lfloor\frac{a_i}{2}\rfloor-1$ and suppose that $G$ is  an $F$-free  graph  of sufficiently large order $n$. If $q(G)\geq q(S_{n,h})$,
 then there exists an induced subgraph $H\subseteq G$ of order $n_1 > n-h^2$ with (i) $\delta(H)\geq h$, (ii)
a vertex $u\in V(H)$ satisfying $d_G(u)=n-1$, and (iii) $d_H(v)\leq h-1$ for every vertex $v\in V(G)\backslash V(H)$.
\end{Lemma}

\begin{Proof}
 If $\delta(G)\geq h$, then let $G=H$. By Lemma~\ref{Lem3.1}, $H$ is the desired graph.
Next we  assume that $\delta(G)\leq h-1$.  We construct a  sequence of graphs $G=G_0\supseteq G_1\supseteq \dots\supseteq G_r$
such that  $\delta(G_i)\leq h-1$ for  $0\leq i\leq r-1$  and $\delta(G_r)\geq h$ as follows.

 Suppose that we have reached such a graph $G_i$ for some $i$ and $\delta(G_i)\leq h-1$.
 Choose $v\in V(G_i)$ with $d_{G_i}(v)=\delta(G_i)$ and let $G_{i+1}=G_i-v$.
Since  $G_i$ is $F$-free and  $|V(G_i)|=n-i$,  Theorem~\ref{edge extremal graph for linear forest} and Lemma~\ref{Lem2.2}~($iii$) imply that
$$h(n-i)-\frac{h^2+h}{2}+1=e(S^+_{n-i,h})\geq e(G_i)\geq e(G)-i(h-1)>h(n-h)-i(h-1).$$
Thus $i<h^2$. Let $H=G_r$, where $r-1$ is the maximum value of $i$ such that $\delta(G_i)\leq h-1$. This implies that $\delta(H)\geq h$.
Moreover, by Lemma~\ref{Lem3.1}, $G$ has a vertex $u$ with $d_G(u)=n-1$.
From the procedure constructing $\{G_i\}^{r}_{i=0}$, we have $u\in V(H)$  and  $d_H(v)\leq h-1$ for  every vertex $v\in V(G)\backslash V(H)$.
\end{Proof}

\begin{Lemma}\label{Lem3.3}
Let   $F=\cup_{i=1}^k P_{a_i}$  be a linear forest with    $k\geq2$,
$a_1\geq \cdots \geq a_k\geq2$, and  at most two odd   $a_i$. Denote $ h=\sum _{i=1}^k \lfloor\frac{a_i}{2}\rfloor-1$ and suppose that $G$ is  an $F$-free  graph  of sufficiently large order $n$.
Let $H$ be an induced subgraph of order $n_1 > n-h^2$ stated in Lemma~$4.2$. 
If $q(G)\geq q(S_{n,h})$,
then $H$ is not a subgraph of $ L_{t_1,t_2,h,h+1}$, where $n_1=t_1h+t_2(h+1)+1$ and $t_1,t_2\geq 0$.
\end{Lemma}

\begin{Proof}
Suppose that $H\subseteq L_{t_1,t_2,h,h+1}$.
Let $H'=G-V(H)$ and $n_2=|V(H')|$.
Let $u$ be a vertex in $H$ with $d_G(u)=n-1$. Since $d_G(u)=n-1$ and $d_H(v)\leq h-1$ for every vertex $v\in V(H')$, we have $d_{H-u}(v)\leq h-2$ for every vertex $v\in V(H')$.
Thus $$e(V(H'),V(H-u))\leq (h-2) n_2< h^3-2h^2.$$
Since $H\subseteq L_{t_1,t_2,h,h+1}$, we have
\begin{eqnarray*}
  e(H) &\leq &e(L_{t_1,t_2,h,h+1})
    =\frac{t_1h^2+t_2(h+1)^2+n_1-1}{2} \\
   &\leq& \frac{t_1h(h+1)+t_2(h+1)^2+n_1-1}{2}
   = \frac{(h+2)(n_1-1)}{2}.
\end{eqnarray*}
Furthermore, since  $q(G)\geq q(S_{n,h})$, by Lemma~\ref{Lem2.2}~($iii$), we have $e(G)>h(n-h)$.
Then
 \begin{eqnarray*}
        e(H')&=& e(G)-e(H)-e(V(H'),V(H-u))-n_2  \\
       & > & h(n-h)-\frac{(h+2)(n_1-1)}{2}-(h^3-2h^2)-n_2\\
        & = &\frac{(h-2)n-(2h^3-2h^2-h-2)}{2}+\frac{h n_2}{2}
         > \frac{h n_2}{2}.
  \end{eqnarray*}
By Theorem~\ref{Thm2.6},  $H'$ contains a path of order $h+2$, denoted by $P$, as a subgraph.
Moreover, since $H\subseteq L_{t_1,t_2,h,h+1}$ and $\delta(H)\geq h$, we have that $H$ is obtained from  $L_{t_1,t_2,h,h+1}$ by deleting some matchings in  its induced
subgraphs $K_{h+2}'s$.
 Hence  $H$ contains $t_1+t_2$ paths of order $h+1$ with end vertex $u$, denoted by $Q_1,\ldots,Q_{t_1+t_2}$, and they have a unique common vertex $u$.
 Note that $n_1$ is sufficiently large, we have $t_1+t_2\geq2$. Since $u$ is the common vertex of $Q_1$ and $Q_2$, we can obtain a longer path $Q$ of order $2h+1$ in $H$ such that $V(Q)=V(Q_1)\cup V(Q_2)$.
Now $G$ contains two disjoint  paths $ Q$   of order $2h+1$ and $P$ as  subgraphs.
We first assume that there are precisely two odd $a_i$, say $a_{i_1}$ and $a_{i_2}$ with $a_{i_1}\geq a_{i_2}$.  Since
  $$h+2=\sum \limits_{i=1}^k \Big\lfloor\frac{a_i}{2}\Big\rfloor+1\geq \frac{a_{i_1}-1}{2}+\frac{a_{i_2}-1}{2}+1\geq a_{i_2},$$
  $$2h+1=2\sum \limits_{i=1}^k \Big\lfloor\frac{a_i}{2}\Big\rfloor-1=\sum \limits_{i=1}^k a_i-3 \geq \sum \limits_{1\leq i\leq k, i\neq i_2}a_i,$$
 we have $F\subseteq G$,  a contradiction.
We next assume that there is precisely one odd  $a_i$.  Note that $d_G(u)=n-1$. We can obtain a path $T$ of order $2h+3$ in $G$¡¡ such that $V(T)=V(P)\cup V(Q_1)$.
Since $$2h+3=2\sum \limits_{i=1}^k \Big\lfloor\frac{a_i}{2}\Big\rfloor+1=\sum \limits_{i=1}^k a_i,$$  we have $F\subseteq G$, a contradiction.
This completes the proof.
\end{Proof}

\begin{Lemma}\label{Lem3.4}
Let $F\in\{P_6\cup P_3, P_7\cup P_3, P_4\cup 2\cdot P_3 \}$ and $G$ be an $F$-free  graph  of  order $n\geq63$.
Let $n_1\geq10$ be even and $H$ be   an induced subgraph of  $G$ of order  $n_1$ with $\delta(H)\geq3$.
 If  $q(G)\geq q(S_{n,3})$, then $H$ is not a subgraph of $ F_{n_1,3}$.
\end{Lemma}
\begin{Proof}
Suppose  that  $H\subseteq F_{n_1,3}$. Since $\delta(H)\geq3$,
 there exists $I\subseteq V(H)$ of size $n_1-2$
 such that $H[I]=\frac{n_1-2}{2}K_2$.  Let $H'=G-V(H)$. If $H'$ contains $P_3$ as a subgraph  or there is an edge with one end vertex in $V(H')$ and
the other in $I$, then   $F\subseteq G$, a contradiction.
So $V(H')\cup I$ induces independent edges and isolated vertices.
It follows that $G\subseteq F_{n,3}$.
By Lemmas~\ref{Lem2.2} ($i$) and \ref{lbound for extremal graph},
 $q(G)\leq q(F_{n,3})\leq n+2<q(S_{n,3})$,  a contradiction.
\end{Proof}

\vspace{2mm}
Although the proof of  Lemma~\ref{Lem3.5} is similar to that of Lemma~\ref{Lem3.4}, we retain it for  completeness.

\begin{Lemma}\label{Lem3.5}
Let $F\in\{P_9\cup P_3, P_6\cup  2\cdot P_3 \}$ and $G$ be an $F$-free  graph  of  order $n\geq112$.
Let $n_1\geq11$ be odd and $H$ be   an induced subgraph  of $G$ of order $n_1$ with $\delta(H)\geq4$.
 If  $q(G)\geq q(S_{n,4})$, then $H$ is  not a subgraph of $F_{n_1,4}$.
\end{Lemma}
\begin{Proof}
Suppose  that  $H\subseteq F_{n_1,4}$. Since $\delta(H)\geq4$,
 there exists $I\subseteq V(H)$ of size $n_1-3$
 such that $H[I]=\frac{n_1-3}{2}K_2$.  Let $H'=G-V(H)$. If $H'$ contains $P_3$ as a subgraph  or there is an edge with one end vertex in $V(H')$ and
the other in $I$, then   $F\subseteq G$, a contradiction.
So $V(H')\cup I$ induces independent  edges and isolated vertices.
It follows that $G\subseteq F_{n,4}$.
By Lemmas~\ref{Lem2.2} ($i$) and \ref{lbound for extremal graph},
 $q(G)\leq q(F_{n,4})\leq n+4<q(S_{n,4})$,  a contradiction.
\end{Proof}

\begin{Lemma}\label{Lem3.6}
Let  $F= P_5\cup P_3$ and $G$ be an $F$-free  graph  of order $n\geq28$. Suppose that $G$ has a vertex $u$ with  $d_G(u)=n-1$.
Let $H$ be an induced subgraph of  $G$ of order $n_1\geq7$  with   $u\in V(H)$ and $\delta(H)\geq2$.
If $q(G)\geq q(S_{n,2})$, then
 $H$ is not a subgraph of $H_{n_1,1}$.
\end{Lemma}

\begin{Proof}
Suppose that $H\subseteq H_{n_1,1}$. Obviously,  $d_H(u)=n_1-1$. Since $\delta(H)\geq2$ and $d_H(u)=n_1-1$,  we have that $H= H_{n_1,1}$.
Let $v$ be a vertex of $H$ with  the second largest degree and  $H'=G-V(H)$.
Note that  $d_G(u)=n-1$.
 If $H'$ contains an edge or there is an edge between  $V(H')$ and $V(H)\backslash \{u,v\}$, then
 $F\subseteq G$, a contradiction. Hence $G\subseteq H_{n,1}$.
 By Lemma~\ref{Lem2.3}~($i$), $q(G)\leq q(H_{n,1})< q(S_{n,2})$,
 a contradiction.
\end{Proof}


\begin{Lemma}\label{Thm3.1}
Let   $F=\cup_{i=1}^k P_{a_i}$  be a linear forest with   $k\geq2$ and
$a_1\geq \cdots \geq a_k\geq2$. Denote $ h=\sum_{i=1}^k \lfloor\frac{a_i}{2}\rfloor-1$ and suppose that $G$ is  an $F$-free  graph  of sufficiently large order $n$.
 If there  exists an even $a_i$  and  at most  two odd  $a_j$, then $q(G)\leq q(S_{n,h})$ with equality if and only if $G=S_{n,h}$.
\end{Lemma}
\begin{Proof}
 Let $G$ be an $F$-free graph of order $n$ with maximum signless spectral radius. It suffices to prove that $G=S_{n,h}$. Since $S_{n,h}$ is $F$-free, we have $q(G)\geq q(S_{n,h})$.  By Lemma~\ref{Lem3.1}, $\Delta(G)=n-1$.   Choose $u\in V(G)$ with $d_G(u)=n-1$.
We first assume that $h=1$. Then  $F=2\cdot P_2$  or  $F=P_2\cup P_3$.
 If there is an edge in $G-u$, then $F\subseteq G$, a contradiction.
This implies that $V(G-u)$ is an independent set and thus  $G=S_{n,1}$.
We next assume that $h\geq2$.

\vspace{2mm}
{\bf Claim.} $\delta(G)\geq h$.

Suppose that $\delta(G)\leq h-1$. By  Lemma~\ref{Lem3.2}, there exists an induced subgraph $H$ of order $n_1>n-h^2$  with  $\delta(H)\geq h$, $u\in V(H)$,
and  $d_H(v)\leq h-1$ for every $v\in V(G)\backslash V(H)$.
By Theorems~\ref{Thm2.1}-\ref{Thm2.4} and Lemmas~\ref{Lem3.3}-\ref{Lem3.5}, $H\subseteq S_{n_1,h}$.
Then there exists $I\subseteq V(H)$ of size $n_1-h$ such that $I$ induces an independent set of $H$
and $u\in V(H)\backslash I$. Since $\delta(H)\geq h$, every vertex in $I$ is adjacent to every vertex in $V(H)\backslash I$. We first assume that there are precisely two odd  $a_i$, say $a_{i_1}$ and $a_{i_2}$ with $a_{i_1}\geq a_{i_2}$.
Note that there exist $3$ disjoint complete bipartite subgraphs $H_i$ of $H$, with $H_i=(X_i, Y_i; E_i)$ for $1\leq i\leq3$,
where $|X_1|=|Y_1|=\sum _{1\leq i\leq k,i\notin\{ i_1,i_2\}}\lfloor \frac{a_i}{2}\rfloor-1$  with $u\in X_1$, $|X_2|=\lfloor\frac{a_{i_1}}{2}\rfloor$, $|Y_2|=\lfloor\frac{a_{i_1}}{2}\rfloor+1$, $|X_3|=\lfloor\frac{a_{i_2}}{2}\rfloor$, and  $|Y_3|=\lfloor\frac{a_{i_2}}{2}\rfloor+1$.
Hence there exist $3$ disjoint  paths $Q_i$ for $1\leq i\leq3$, where $Q_1$ is path of order  $\sum _{1\leq i\leq k,i\notin\{ i_1,i_2\}} a_i-2$  in $H_1$ with one end vertex $u$,
 $Q_2$ is a path of order  $a_{i_1}$   in $H_2$, and $Q_3$ is a path of order  $a_{i_2}$   in $H_3$.
Let $H'=G-V(H)$.
If $V(H')\cup I$ induces at least an edge $u_1u_2$ (if one of $u_1$ and $u_2$ is in $Y_i$, say $u_1\in Y_i$, then we can find another $u_1'\in I\backslash \cup_{i=1}^3Y_i$ to substitute $u_1$ in $Y_i$), then  the path $Q_1$ of order $\sum _{1\leq i\leq k,i\notin\{ i_1,i_2\}} a_i-2$ can be extended to a path $Q_1u_1u_2$ of order
$\sum _{1\leq i\leq k,i\notin\{ i_1,i_2\}} a_i$ in $G$. Hence $F\subseteq G$, a contradiction.
This implies that $V(H')\cup I$ is an independent set and thus   $G\subseteq S_{n,h}$. Next assume that there is precisely one odd $a_i$. By similar arguments , we can get the same result that $V(H')\cup I$ is an independent set and thus   $G\subseteq S_{n,h}$.  Since $\delta(G)\leq h-1$,  $G$ is a  proper subgraph of $ S_{n,h}$. By Perron--Frobenius  theorem, $q(G)<q(S_{n,h})$, a contradiction. This completes the claim.

By Claim, $\delta(G)\geq h$. It follows from Theorems~\ref{Thm2.1}-\ref{Thm2.4}, Lemmas~\ref{Lem2.2}, and Lemma~\ref{lbound for extremal graph} that $G\subseteq S_{n,h}$. By Perron--Frobenius  theorem and the extremality of $G$, $G=S_{n,h}$.
\end{Proof}

\begin{Lemma}\label{Thm3.2}
Let $F= \cup_{i=1}^2 P_{a_i}$ with  $a_1\geq a_2$ and $a_1>3$. Denote $h=\sum_{i=1}^2 \lfloor\frac{a_i}{2}\rfloor-1$
 and suppose that $G$ is an $F$-free  graph  of  sufficiently large order  $n$.
 If $a_1$ and $a_2$ are odd,  then $q(G)\leq q(S^+_{n,h})$ with equality if and only if $G= S^+_{n,h}$.
\end{Lemma}

\begin{Proof}
 Let $G$ be an $F$-free graph of order $n$ with maximum signless spectral radius. It suffices to prove that $G=S^+_{n,h}$. Since $S^+_{n,h}$ is $F$-free, we have $q(G)\geq q(S^+_{n,h})$. Obviously, $q(G)\geq q(^+_{n,h})> q(S_{n,h})$.
 By  Lemma~\ref{Lem3.1}, $\Delta(G)=n-1$. Choose $u\in V(G)$ with $d_G(u)=n-1$.

\vspace{2mm}
{\bf Claim.} $\delta(G)\geq h$.

 Suppose that $\delta(G)\leq h-1$. By Lemma~\ref{Lem3.2}, there exists an induced subgraph $H$ of order $n_1>n-h^2$ with $\delta(H)\geq h$, $u\in V(H)$, and $d_H(v)\leq h-1$ for every $v\in V(G)\backslash V(H)$.
By Theorems~\ref{Thm2.3}-\ref{Thm2.4} and  Lemmas~\ref{Lem3.3}-\ref{Lem3.6},
$H\subseteq S^+_{n_1,h}$. Thus there exists $I\subseteq V(H)$ of size $n_1-h$  such that $I$ induces at most one edge of $H$
and  $u\in V(H)\backslash I$.   Let $H'=G-V(H)$.

{\bf Case 1.}  $I$ induces precisely one edge  $v_1v_2$ and isolated vertices.
  Since $\delta(H)\geq h$, every vertex in $I\backslash\{v_1,v_2\}$ is adjacent to every vertex in $V(H)\backslash I$ and
 $v_2$ is adjacent to at least $h-1\geq \lfloor\frac{a_1}{2}\rfloor-1$ vertices in  $V(H)\backslash I$.
Then there exist two disjoint complete bipartite subgraphs $H_i$ of $H$, with $H_i=(X_i, Y_i; E_i)$ for $1\leq i\leq 2$,
where $|X_1|=|Y_1|=\lfloor\frac{a_1}{2}\rfloor-1$ with $u\in X_1$, $v_2\in Y_1$,  $v_1\notin Y_1\cup Y_2$, $|X_2|=\lfloor\frac{a_2}{2}\rfloor$, and  $|Y_2|=\lfloor\frac{a_2}{2}\rfloor+1$.
Hence there exist two disjoint  paths $Q_i$ for $1\leq i\leq 2$, where $Q_1$ is path of order  $a_1-3$ with end vertices $u$ and $v_2$  in $H_1$ and
 $Q_2$ is a path of order  $a_2$   in $H_2$.
If $V(H')$ induces at least an edge $w_1w_2$, then  the path $Q_1$ of order $a_1-3$ with end vertices $u$ and $v_2$  can be extended to a path $v_1Q_1w_1w_2$ of order
$a_1$ in $G$.  Hence $F\subseteq G$, a contradiction. This implies that  $V(H')$ induces an independent set.
Similarly,  we can prove that there is no edge with one end vertex in $I$ and the other  in $V(H')$.
Hence  $V(H')\cup I$ induces precisely one edge $v_1v_2$ and isolated vertices.
Moreover, since $\delta(G)\leq h-1$,  $G$ is a proper subgraph of   $S^+_{n,h}$.  By Perron--Frobenius  theorem, $q(G)<q(S^+_{n,h})$,   a contradiction.

{\bf Case 2.}  $I$ is an independent set.
 Since $\delta(H)\geq h$,   every vertex in $I$
is adjacent to every vertex in $V(H)\backslash I$. If $V(H')\cup I$ induces a path of order 3, then by similar arguments to the case above,
we can find two disjoint paths of order $a_1$ and $a_2$ respectively. Then $F \subseteq G$, a contradiction.
Hence $V(H')\cup I$ induces  independent edges  and isolated vertices, i.e., $v_1v_2,\ldots, v_{2p-1}v_{2p}, w_1,\ldots, w_q$.
Since $\delta(G)\leq h-1$, if  $p\leq1$ then
  $G$ is a proper subgraph of $ S^+_{n,h}$ and thus
 $q(G)<q(S^+_{n,h})$, a contradiction.
So  $p\geq2$. Then all $w_i$ are not adjacent to any vertex in $I$, i.e., all $v_{2i-1}v_{2i}\in E(H')$, otherwise
 we can find two disjoint paths of order $a_1$ and $a_2$ respectively, and thus $F \subseteq G$.
 Likewise, we see that each vertex of all edges of $H'$ is adjacent to none of vertices in $V(H)\backslash I$ except for $u$.
 We show  that $q(G)<q(S^+_{n,h})$.
 Let $\mathbf{x}=(x_1,x_2,\dots,x_n)$ be a positive unit eigenvector of $Q(G)$ corresponding to $q(G)$. Choose a vertex $v\in V(H)\backslash (I\cup \{u\})$. 
For simplicity, let $q=q(G)$. From the eigenequations  of $Q(G)$, we have  for $1\leq i\leq p$,
 $$(q-2)x_{v_{2i-1}}=x_{v_{2i}}+x_u, \quad(q-2)x_{v_{2i}}=x_{v_{2i-1}}+x_u,$$
 $$(q-d(v))x_v=\sum \limits_{w\in N(v)}x_w>x_u.$$
Since $d(v)> n-h^2-h>3$,  we have $$x_{v_{2i-1}}=x_{v_{2i}}=\frac{x_u}{q-3}, \qquad x_v>\frac{x_u}{q-d(v)}>x_{v_{2i-1}}.$$
Let $G'$ be a graph obtained from $G$ by deleting all edges in $\{v_{2i-1}v_{2i}: 1\leq i\leq p\}$ and
adding all edges in  $\{vv_i: 1\leq i\leq 2p\}$.  Obviously, $G'\subseteq S_{n,h}$ and
$$q(G')\leq q(S_{n,h})<q(S^+_{n,h}).$$ On the other hand,
$$q(G')-q(G)=\sum\limits_{1\leq i \leq2p}(x_v+x_{v_i})^2-\sum\limits_{1\leq i\leq p}(x_{v_{2i-1}}+x_{v_{2i}})^2
>\sum\limits_{1\leq i \leq p}(x_v+x_{v_{2i-1}})^2>0,$$
implying that $q(G)<q(G')$. Thus $q(G)<q(S^+_{n,h})$, a contradiction.
This completes the claim.
\vspace{2mm}

By Claim, $\delta(G)\geq h$. It follows from Theorems~\ref{Thm2.3}-\ref{Thm2.4} and Lemmas~\ref{Lem2.2}-\ref{lbound for extremal graph} that $G\subseteq S^+_{n,h}$.  By Perron--Frobenius  theorem and the extremality of $G$, $G=S^+_{n,h}$.
\end{Proof}

\begin{Lemma}\label{spectral radius for kP3}
Let $k\geq2 $ and $G$ be a $k\cdot P_3$-free graph of order $n\geq \frac{11}{2}k^2-9k+12$. Then $q(G)\leq q(F_{n,k})$ with equality if and
only if $G=F_{n,k}$.
\end{Lemma}
\begin{Proof}
%
Let $G$ be a $k\cdot P_3$-free graph of order $n$ with maximum signless spectral radius. It suffices to prove that $G=F_{n,k}$. Since $F_{n,k}$ is  $k\cdot P_3$-free, we have $q(G)\geq q(F_{n,k})$. By Lemma~\ref{lbound for extremal graph}, $q(G)>n+2k-5$.

First we assume that $k=2$. By Lemma~\ref{2P3-free}, $G\subseteq F_{n,1}(K_h)$ with $4\leq h\leq5$, $G\subseteq F_{n,1}(N_6)$, or  $G\subseteq F_{n,2}$.
By Perron--Frobenius  theorem, Lemma~\ref{Kh}, and Lemma~\ref{N1}, $G=F_{n,2}$.
So next we assume that $k\geq3$.
Note that
$$q(G)\leq \max \limits_{v\in V(G)}\Big\{d(v)+\frac{1}{d(v)}\sum\limits_{z\in N(v)} d(z)\Big\},$$
which dates back to Merris \cite{ Merris}.
Let $u$ be a vertex of $G$ such that
$$q(G)\leq \max \limits_{v\in V(G)}\bigg\{d(v)+\frac{1}{d(v)}\sum\limits_{z\in N(v)} d(z)\bigg\}=d(u)+\frac{1}{d(u)}\sum\limits_{z\in N(u)} d(z).$$
 We first claim that $d(u)\geq 2k-3$. Otherwise
$$q(G)\leq d(u)+\Delta(G)\leq 2k-4+n-1=n+2k-5,$$ a contradiction.
Note that $q(G)>n+2k-5$, we have
$$\sum\limits_{v\in N(u)} d(v)>(n+2k-5-d(u))d(u).$$
Next we consider two cases.

\vspace{2mm}
{\bf Case 1.} $d(u)\leq \frac{n+k-5}{4}$.  Let $W=V(G)\backslash (\{u\}\cup N(u))$. Let $C\subseteq N(u)$ be the vertex subset such that every vertex in $C$
has at least $2k$ neighbours in $W$.  We claim that $|C|\geq k$. Otherwise, we have
\begin{eqnarray*}
  \sum\limits_{v\in N(u)} d(v) &=& \sum\limits_{v\in C} d(v)+\sum\limits_{v\in N(u)\backslash C} d(v) \\
   &=& |C|(n-1)+(d(u)-|C|)(2k-1+d(u)) \\
   &=& (n-d(u)-2k)|C|+(2k+d(u)-1)d(u)\\
   &\leq& (n-d(u)-2k)(k-1)+(2k+d(u)-1)d(u) \\
   &=& (d(u)(d(u)+k)+(n-2k)(k-1) \\
   &=& (n+2k-5-d(u))d(u)+d(u)(2d(u)-n-k+5)+(n-2k)(k-1)\\
   &\leq& (n+2k-5-d(u))d(u)+(2k-3)(-n+3k-1)+(n-2k)(k-1)\\
    &=&(n+2k-5-d(u))d(u)-((k-2)n-4k^2+9k-3)\\
   & \leq &(n+2k-5-d(u))d(u),
\end{eqnarray*}
where the last third inequality holds  because the function $f(x)=x(2x-n-k+5)$ is decreasing with respect to $x$ for
$2k-3\leq x\leq \frac{n+k-5}{4}$, the last inequality holds  as  $n\geq \frac{11}{2}k^2-9k+12$,
a contradiction.  Then we can embed $k\cdot P_3$ in $G$  with all centers in $C$ and all ends in $W$,  a contradiction.

\vspace{2mm}
{\bf Case 2.} $d(u)>\frac{n+k-5}{4}$. Let $G_u=G[V(G)\backslash \{u\})]$.
Then $G_u$ must be $(k-1)\cdot P_3$-free. Otherwise, $(k-1)\cdot P_3$ in $V(G)\backslash \{u\})$ and  a disjoint $P_3$ with center $u$ and two ends in $N(u)$  yield  $k\cdot P_3$, a contradiction.
%
In addition, by  inequalities $q(G) \leq \frac{2e(G)}{n-1}+n-2$ \cite{FY,YF} and  $q(G)>n+2k-5$,
we have  $e(G)>(k-\frac{3}{2})(n-1)$. Then
$$e(G_u) = e(G)-d(u)  > (k-\frac{3}{2})(n-1)-(n-1)=(k-\frac{5}{2})(n-1).$$
Note that $G_u$ is $(k-1)\cdot P_3$-free.
By Theorem~\ref{stablility thm for kP3}, we have (i) $G_u\subseteq F_{n-1,k-1}$,
(ii) $G_u\subseteq F_{n-1,k-2}(K_h)$ for $4\leq h\leq5$, or (iii) $G_u\subseteq F_{n-1,k-2}(N_6)$.
In (i), 
we have $G\subseteq F_{n,k}$.  By Perron--Frobenius  theorem and the extremality of $G$, we have $G=F_{n,k}$.
In (ii), we have $G\subseteq F_{n,k-1}(K_h)$ for
 $4\leq h\leq5$. By Lemma~\ref{Kh}, $$q(G)\leq q(F_{n,k-1}(K_h))<q(F_{n,k}),$$ a contradiction.
 In (iii), we have $G\subseteq F_{n,k-1}(N_6)$.
 By Lemma~\ref{N1}, $$q(G)\leq q(F_{n,k-1}(N_6))<q(F_{n,k}),$$ a contradiction. This completes the proof.
\end{Proof}

\vspace{2mm}
Now we are ready to prove Theorem~\ref{Q spectral radius for all linear forest}.

\vspace{2mm}
\noindent {\bf Proof of Theorem~\ref{Q spectral radius for all linear forest}.} Theorem \ref{Q spectral radius for all linear forest} directly follows from
Lemmas~\ref{Thm3.1}-\ref{spectral radius for kP3}.\QEDB

\vspace{3mm}

\noindent{\bf Acknowledgements.} The authors would like to thank the anonymous referees for many helpful  suggestions to an earlier version of this paper, which results in a great improvement.

\end{document}